\documentclass[11pt,a4paper]{article}
\usepackage{amssymb}
\usepackage{subfloat,amsmath,graphics,amssymb,amsfonts,graphicx,lineno}
\usepackage{amsmath}
\usepackage{amssymb}
\usepackage{amsthm}
\usepackage{amsfonts}
\newtheorem{theorem}{Theorem}[section]

\newtheorem{lemma}[theorem]{Lemma}

\theoremstyle{definition}

\theoremstyle{definition}

\begin{document}
\begin{center}

\textbf{Perturbations of Banach algebras and amenability}
\\ \ \\
Miad Makareh Shireh
\ \\ \ \\
\end{center}
ABSTRACT: In this paper we prove that if $(A,\pi)$ is an amenable Banach algebra and if  $\rho$ is another Banach algebra multiplication on $A$ such that $\|\rho-\pi\|<{1 \over 11}$, then $(A, \rho)$ is also amenable.

 \section{Introductions}
 Let $A$ to be a Banach algebra and $X$ an $A$-bimodule that is a Banach space. We say that $X$ is a Banach $A$-bimodule if there exists  constant $C>0$ such that
 \begin{align*}
 \|a.x\|&\leq C\|a\|\|x\|,\\
 \|x.a\|&\leq C\|a\|\|x\| \hskip 1 cm (a \in A, x \in X).
 \end{align*}
 If $X$ is a Banach $A$-bimodule, then $X^*$ is a Banach $A$-bimodule for the actions defined by \begin{align*}
\langle a.f,x\rangle &=\langle f,x.a\rangle\\ \langle f.a,x\rangle &=\langle f,a.x\rangle \hskip 0.5cm (a\in A, f\in X^*,x \in X). \end{align*}
The Banach $A$-bimodule $X^*$ defined in this way is said to be a dual Banach $A$-bimodule.\\\\

 A linear mapping $D$ from $A$ into $X$ is a derivation if $$ D(ab)=a.D(b)+D(a).b \hskip 1 cm (a,b \in A).$$
For $x \in X$, the mapping $ad_x:A\longrightarrow X$ defined by $ad_x(a)=a.x-x.a$ is a continuous derivation. The derivation $D$ is inner if there exists $x \in X$ such that $D=ad_x$.\\
 $A$ is said to be amenable if for every Banach $A$-bimodule $X$ , any continuous derivation from $A$ into the dual Banach $A$-bimodule $X^*$ is inner. This notion has been introduced in [2] and has been studied extensively.\\\\
 Let $A$ be an Banach algebra. $A^{op}$ is another Banach algebra which is the same as $A$ as Banach spaces but the product of $A^{op}$ is the reverse of the product of $A$ i.e. $$a\circ b=ba \hskip 1 cm (a,b \in A),$$
 where $\circ$ denotes the multiplication of $A^{op}$.

 The so-called multiplication map, denoted by $\pi$,  $\pi:A\hat{\otimes}A^{op}\longrightarrow A$  is specified by  $$\pi(a\otimes b)=ab \hskip 1 cm (a,b\in A)$$

 By the difference between the two multiplications $\pi$ and $\rho$ on a Banach algebra $A$, we mean the norm of $\pi-\rho$ as an operator from $A\hat{\otimes}A^{op}$ to $A$.
 In [3] Johnson proved that if $(A,\pi)$ is amenable, then there exists an $\epsilon>0$ such that if $\rho$ is another Banach algebra multiplication on  $A$ such that on   $\|\pi-\rho\|<\epsilon$ , then $(A, \rho)$ is also amenable. But that $\epsilon$ here depends on the structure of the Banach algebra $A$. In this paper we give a partially different proof for that theorem and we prove the following result:\\\\
 If $(A,\pi)$ is an amenable Banach algebra, then $(A,\rho)$ is also amenable for every Banach algebra multiplication $\rho$ on $A$ such that $\| \pi-\rho\|<{1 \over 11}$.

\section{Perturbations of Banach algebras}
Before going to the mail theorem, we bring two lemmas from [3] that are used in our proof. \\\\
For two closed subspaces $Y$ and $Z$ of a Banach space $X$,their Hausdorff distance is defined by $$d(Y,Z)= \max\{\sup \{d(y,Z):\|y\|\leq 1\},\sup\{d(z,Y): \|z\|\leq 1\}\}$$
\begin{lemma}
Let $Y$ and $Z$ be closed subspaces of a Banach space $X$. Suppose that there is a projection $P$ of $X$ onto $Y$ with $\|P\|< d(Y,Z)^{-1}-1$. Then $P$ maps $Z$ one to one onto $Y$ and the inverse $\alpha$ of $P|_Z$satisfies $(d=d(Y,Z))$
\begin{align*}
\|\alpha \| &\leq (1+d)(1-\|P\|d)^{-1}\\
\| \alpha(y)-y \| &\leq ((1+d)(1-\|P\|d)^{-1}-1)\|y\|\\
\|P(z)-z\| &\leq d(1+\|P\|)\|z\|
\end{align*}

\end{lemma}
{\bf Proof:} See [3, Lemma 5.2].
\begin{lemma}
Let $X_1$ and $X_2$ be Banach spaces and $S,T \in B(X_1,X_2)$ and let $S$ be onto. Suppose that there exists $K>0$ such that for all $y\in X_2$, there is $x \in X_1$ with $\|x\|\leq K\|y\|$ and $S(x)=y$. If $K\|S-T\|<1$, then $T$ will also be onto and for each $y \in X_2$, there exists $x \in X_1$ such that $\|x\|\leq K(1-K\epsilon)^{-1}\|y\|$ and $T(x)=y$, where $\epsilon=\|S-T\|$.
\end{lemma}
{\bf Proof:} It is a special case of [3, Lemma 6.1].\\\\
In next theorem and note we denote all  multiplications induced by $\pi$ by a sign of $\pi$ for example in order to show the product of $a$ and $b$ induced by $\pi$, we use $a_{\pi}b$, We have the same way to show them for $\rho$.
{\bf Note:} If $\pi^{\#}$ and $\rho^{\#}$ are the products respectively induced by $\pi$ and $\rho$ on $A^{\#}$ ($A^{\#}$ is the unitization of $A$) then we have $$\|(\pi^{\#}-\rho^{\#})((a,\alpha)\otimes(b,\beta))\|=\|a_{\pi}b-a_{\rho}b\|\leq \|\pi-\rho\|\|a\|\|b\| \hskip 0.5 cm (a,b \in A).$$
And hence $$\|(\pi^{\#}-\rho^{\#})((a,\alpha)\otimes(b,\beta))\|\leq \|\pi-\rho\|\|(a,\alpha)\|\|(b,\beta)\|$$
Thus we have $$\|\pi^{\#}-\rho^{\#}\|\leq \|\pi-\rho \|.$$

\begin{theorem}
Let $(A,\pi)$ be an amenable Banach algebra. If $\rho$ is another Banach algebra multiplication on $A$ such that $\|\pi- \rho\|<{1 \over 11}$, then $(A,\rho)$ is also amenable.
\end{theorem}

{\bf Proof:}{\it By} the note above, we can assume that $A$ has and identity $1$ for both multiplications $\pi$ and $\rho$. Let  $j:A\longrightarrow A\widehat{\bigotimes}A$ be defined by $j(a)=a\otimes1$.\\
Then $\|j\|\leq1$ and $\pi j=Id_A$. So $\pi^{**}j^{**}=Id_{A^{**}}$ . It can be easily checked that $P={\rm Id}_{(A\widehat{\bigotimes}A)^{**}}-j^{**}\pi^{**}$ is a projection onto ker$\pi^{**}$ with norm at most 2.\\
By Lemma 2.2, and letting $X_1=(A\widehat{\bigotimes}A)^{**}$and $X_2=A^{**}$, $S_1=\pi^{**}$,$T_1=\rho^{**}$,by $K=1$ $(\rm{since} \|j^{**}\|\leq1)$, we get that for $\|S_1-T_1\|=\epsilon<1$, $\rho^{**}$ will be onto and for every $F \in \rm{ker}\pi^{**}$, there is $B \in (A\widehat{\bigotimes}A)^{**}$ such that $\rho^{**}(B)=\rho^{**}(F)$ and $$ \|B\|\leq (1-\epsilon)^{-1}\|\rho^{**}(F)\|=(1-\epsilon)^{-1}\|\rho^{**}(F)-\pi^{**}(F)\|\leq (1-\epsilon)^{-1}\epsilon\|F \|$$
So $F-B \in \rm{ker} \rho^{**}$ and $\|F-(F-B)\|=\|B\|\leq \epsilon(1-\epsilon)^{-1}\|F\| $. So that
$$\sup\{d(F,\rm{ker}\rho^{**}) : F \in \rm{ker} \pi^{**} and \|F\|\leq1\}\leq \epsilon(1-\epsilon)^{-1}.$$
And similarly by changing the role of $S_1$ and $T_1$, we will obtain
$$\sup\{d(F,ker\pi^{**}) :F \in \rm{ker} \rho^{**} and \|F\|\leq1\}\leq\epsilon(1-\epsilon)^{-1}$$
Hence $$d:=d(\rm{ker} \pi^{**},\rm{ker} \rho^{**}) \leq \epsilon(1-\epsilon)^{-1}.$$
So if $\epsilon<{1 \over 4}$, then $$\|P\|\leq 2<(\epsilon(1-\epsilon)^{-1})^{-1}-1\leq d(\rm{ker} \pi^{**},\rm{ker} \rho^{**})^{-1}-1.$$
And hence by Lemma 2.1, there exists a linear homeomorphism $\alpha$ from $\rm{ker} \pi^{**}$ onto $\rm{ker} \rho^{**}$ such that $$\|\alpha\|\leq (1-3\epsilon)^{-1}, \|\alpha^{-1}\| \leq\|P\|\leq2$$
$$\|F-\alpha(F)\|\leq 3\epsilon(1-3\epsilon)^{-1}\|F \| \hskip 0.5 cm (F \in \rm{ker} \pi^{**})$$
$$\|F-\alpha^{-1}(F)\|\leq 3\epsilon(1-\epsilon)^{-1} \|F\| \hskip 0.5 cm (F \in \rm{ker} \rho^{**}).$$
Suppose that $F \in (A\widehat{\otimes}A)$ is an elementary tensor say $b\otimes c$ for $b,c \in A$. Then for $a \in A$, we have \begin{align*}
\|a._{\pi}F-a._{\rho}F\| &=\|a.(b\otimes c)-a._{\rho}(b\otimes c)\|\\ &=\|ab\otimes c-a_{\rho} b\otimes c\|=\|(a_{\rho} b-ab)\| \|c\|\\ &\leq \|\rho-\pi\|\|a\otimes b\|\|c\|\\ &\leq \epsilon \|a\|\|b\|\|c\|=\epsilon\|a\|\|F\|.
\end{align*}
So that $$\|a._{\pi}F-a._{\rho}F\|\leq \epsilon\|a\| \|F\| \hskip 0.5 cm (a \in A,F \in A\widehat{\otimes}A).$$
And by using Goldsteine's Theorem, we have $$\|a._{\pi}F-a._{\rho}F \| \leq \epsilon \|F \| \hskip 0.5 cm (F \in (A\widehat{\otimes}A)^{**}) \hskip 3 cm (\dag) $$
Similarly
$$\|F._{\pi}a-F.{\rho}a\|\leq \epsilon \|a\|\|F\| \hskip 0.5 cm (a \in A, F \in (A\widehat{\otimes}A)^{**}).$$
Now consider the derivation $D:A\longrightarrow \rm{ker} \pi^{**}(\cong (\rm{ker}\pi)^{**})$ by $D(a)=a\otimes1-1\otimes a$, then amenability of $(A,\pi)$ implies the existence of an element $\xi \in \rm{ker}\pi^{**}$ such that $$a\otimes1-1\otimes a=a._{\pi} \xi-\xi._{\pi} a \hskip 0.5 cm (a \in A).$$
Let $\delta=\alpha(\xi) \in \rm{ker} \rho^{**}$. Then we have
\begin{align*}
\|a._{\pi}\xi-a._{\rho}\delta\| &=\|a._{\pi}\xi-a._{\rho}(\alpha(\xi))\| \\ &\leq\|a._{\pi}\xi-a._{\pi}(\alpha(\xi))\|+\|a._{\pi}(\alpha(\xi))-a._{\rho}(\alpha(\xi))\|\\ &\leq 3\epsilon(1-3\epsilon)^{-1}\|a\|\|\xi\|+\epsilon(1-3\epsilon)^{-1}\|a\|\|\xi\|. \hskip 0.3 cm (\rm{By} \hskip 0.1 cm  \rm{properties} \hskip 0.1 cm of \hskip 0.1 cm\alpha \hskip 0.1 cm  \rm{and} \hskip 0.1 cm (\dag))
\end{align*}
And similarly
$$\|\xi._{\pi} a -\delta._{\rho} a \|\leq 4\epsilon(1-3\epsilon)^{-1}\|a\|\|\xi\|.$$
So that
\begin{align*}
\|a\otimes1-1\otimes a-(a._{\rho}\delta-\delta._{\rho}a)\| &=\|a._{\pi} \xi-\xi._{\pi} a-(a._{\rho}\delta-\delta._{\rho}a)\|\\
&\leq \| a._{\pi}\xi-a._{\rho}\delta\|+\|\xi._{\pi} a -\delta._{\rho} a \| \\ &\leq 8\epsilon(1-3\epsilon)^{-1}\|a\|.
\end{align*}
So $$\|a\otimes1-1\otimes a-(a.{\rho}\delta-\delta._{\rho}a)\|\leq O(\epsilon)\|a\| \hskip 1 cm (a \in A). \hskip 1 cm (\ddag)$$
Where $O(\epsilon)\longrightarrow 0$ as $\epsilon\longrightarrow0^+$.\\
From now on all the multiplications we consider are respect to the multiplication $\rho$ on $A$.
We denote the multiplication in $A\widehat{\otimes}A^{op}$ by $\star_{\rho}$. Also we show the Arens product on $(A\widehat{\otimes}A^{op})^{**}$ with the same notation. So for elementary tensors, $$(a\otimes b)\star_{\rho}(c\otimes d)=ac\otimes db$$
 For $R=\sum_i a_i\otimes b_i \in \rm{ker} \rho$ we have
\begin{align*}
R\star_{\rho}\delta-R &= \sum_i (a_i\otimes b_i)\star_{\rho} \delta-\delta\sum_ia_ib_i-\sum_i a_i\otimes b_i +1\otimes\sum_i a_ib_i\\
&=\sum_i (a_i._{\rho}\delta-\delta._{\rho} a_i-a_i\otimes1+1\otimes a_i)._{\rho}b_i.
\end{align*}
So
\begin{align*}
\|R\star_{\rho}\delta-R\| &=\|\sum_i (a_i._{\rho}\delta-\delta._{\rho} a_i-a_i\otimes1+1\otimes a_i)._{\rho}b_i\|\\
&\leq \sum_i\|{a_i \over \|a_i\|}._{\rho}\delta-\delta._{\rho}{a_i\over \|a_i\|}+{a_i \over \|a_i\|}\otimes1+1\otimes{a_i \over \|a_i\|}\|\|a_i\|\|b_i\| \\ &\leq \|R\|\sup_{a \in A_1} \|a._{\rho}\delta-\delta._{\rho} a-a\otimes1+1\otimes a\|.
\end{align*}
Now if $R \in (\rm{ker} \rho)^{**}$, then by Goldsteine's Theorem, there exists a net $(r_i)_i$ with $\|r_i\|\leq \|R\|$, in $\rm{ker} \pi$ such that $r_i\longrightarrow_i R \hskip 0.5 cm {\rm wk^*}$.
Note that since $\rm{ker} \rho^{**}\cong (\rm{ker} \rho)^{**}$,isometrically, then for notational convenience, we don't disguise between $\delta$ as an element in $\rm{ker} \rho^{**}$ and its image as an element of $(\rm{ker} \rho)^{**}$. \\
Thus $$r_i._{\rho}\delta-r_i\longrightarrow_i R._{\rho}\delta-R \hskip 0.5 cm {\rm wk^*}.$$
And hence $\| R._{\rho}\delta-R\|\leq \sup_i \|r_i._{\rho}\delta-r_i\|$. So we have
$$\|R\star_{\rho}\delta-R\|\leq \|R\|\sup_{a \in A_1} \|a._{\rho}\delta-\delta._{\rho} a-a\otimes1+1\otimes a\| \hskip 0.3 cm (R \in (\rm{ker} \rho)^{**}).$$
And hence by $(\ddag)$, we obtain $$\|R\star_{\rho}\delta-R\|\leq O(\epsilon)\|R\| \hskip 0.5 cm (R \in (\rm{ker} \rho)^{**}).$$
If we define $\lambda:(\rm{ker} \rho)^{**}\longrightarrow (\rm{ker} \rho)^{**}$ by $\lambda(S)=S \star_{\rho}\delta$, then for $\epsilon<{1\over11}$, \\ $O(\epsilon)={ 8\epsilon \over (1-3\epsilon)}<1$ and hence $\|\lambda-Id_{(\rm{ker} \rho)^{**}}\|<1$ and thus $\lambda$ will be invertible. \\
 Since $\lambda$ is surjective,  there exists $ x \in (\rm{ker} \rho)^{**}$ such that $\lambda(x)=\delta$. So  $x \star_{\rho}\delta=\delta$ and therefore for every $ y \in (\rm{ker} \rho)^{**}$, we have $(y\star_{\rho}x-y) \star_{\rho}\delta =0$ but this means that
$$\lambda(y \star_{\rho}x-y)=0 \hskip 1 cm (y \in(\rm{ker} \rho)^{**}).$$
Now by injectivity of $\lambda$, we have
 $$ y \star_{\rho}x=y \hskip 1 cm (y \in (\rm{ker} \rho)^{**}).$$
Hence $x$ will be a right identity for $(\rm{ker} \rho)^{**}$ and hence $\rm{ker} \rho$ has a bounded right approximate identity. So from [1,Theorem 3.10] , $(A,\rho)$ is amenable. \hskip 1 cm $\Box$


\begin{thebibliography}{999}
\bibitem[1]{a}P.C.Curtis, R.J.Loy,The structure of amenable Banach algebras, Journal of London Mathematical Society, 40(2)(1989) 89-104.
\bibitem[2]{a}B.E.Johnson, Cohomology in Banach Algebras.  American Mathematical Society, Providence, RI, (1972).
\bibitem[3]{a}B.E.Johnson, Preturbations of Banach algebras. Proc. London Math. Soc. (3) 34 (1977), no 3, 439-458.
\end{thebibliography}
  \end{document}